\documentclass{amsart}
\usepackage{latexsym,amsmath,amssymb}

\theoremstyle{definition}

\theoremstyle{remark}

\numberwithin{equation}{section}

\begin{document}

\title[WEIGHTED CONDITIONAL EXPECTATION]
{UNBOUNDED WEIGHTED CONDITIONAL EXPECTATION OPERATORS}

\author{\sc y. estaremi}
\address{\sc  y. estaremi}
\email{estaremi@gmail.com}
\address{Department of mathematical
sciences, Payame noor university, Tehran, iran}

\thanks{}

\thanks{}



\date{}

\dedicatory{}

\commby{}

\begin{abstract}
In this note basic properties of unbounded weighted conditional expectation operators are investigated. A description of polar decomposition and quasinormality in this context are provided. Also, we give some necessary and sufficient conditions for $EM_u$ to leave invariant it's domain. Finally,  some examples are provided to illustrate concrete application of the main results of the paper.
\noindent {}
\end{abstract}

\maketitle

\section{ \sc\bf Introduction and Preliminaries}
Theory of weighted conditional expectation type operators is one of important
arguments in the connection of operator theory and measure theory.
Weighted conditional expectation type operators have been studied in an operator
theoretic setting, by, for example, De pagter and Grobler \cite{g}
and Rao \cite{rao1, rao2}, as positive operators acting on
$L^p$-spaces or Banach function spaces. In \cite{mo}, S.-T. C. Moy
characterized all operators on $L^p$ of the form $f\rightarrow
E(fg)$ for $g$ in $L^q$ with $E(|g|)$ bounded.  Also, some results
about weighted conditional expectation type operators can be
found in \cite{dou, her, lam}. In \cite{dhd} P.G. Dodds, C.B.
Huijsmans and B. de Pagter showed that lots of operators are of
the form of weighted conditional type operators. The class of
weighted conditional type operators includes the classes of composition
operators, multiplication operators, weighted composition
operators, some integral type operators and etc. These are some
reasons that stimulate us to consider weighted conditional type
operators in our work. So, in \cite{e,e1,e2,e3,ej, ej2} we studied bounded
multiplication conditional expectation
operators $M_wEM_u$ on $L^p$ spaces.\\

 As far as I know, unbounded weighted conditional type operators weren't investigated till now. This paper is devoted to the study of unbounded weighted conditional type operators on $L^2(X,\Sigma,\mu)$, such that $(X,\Sigma,\mu)$ is a $\sigma$-finite measure space. In general closed operators are important classes of unbounded linear
operators which are large enough to cover all interesting operators occurring
in applications. In section 2, we develop basic concepts and results
about weighted conditional type operators on Hilbert space $L^2(\Sigma)$. First we investigate some basic properties of unbounded weighted conditional type operators. We show that a densely-defined weighted conditional type operators are always closed. In general closed operators are important classes of unbounded linear operators which are large enough to cover all interesting operators occurring
in applications.\\

Let $\mathcal{H}$ be stand for a Hilbert space and $\mathcal{B}(\mathcal{H})$ for the Banach algebra of all bounded operators
on $\mathcal{H}$. By an operator in $\mathcal{H}$ we understand a
linear mapping $T:\mathcal{D}(T)\subseteq \mathcal{H}\rightarrow
\mathcal{H}$ defined on a linear subspace $\mathcal{D}(T)$ of
$\mathcal{H}$ which is called the domain of $T$. Set
$\mathcal{D}^{\infty}(T)=\cap^{\infty}_{n=1}\mathcal{D}(T^n)$. Denote by $\mathcal{N}(T)$, $\mathcal{R}(T)$ and $T^{\ast}$ the kernel, the range and the adjoint of $T$ respectively. Given an operator $T$ in $\mathcal{H}$, we write the graph norm
$\|.\|_{T}$ on $\mathcal{D}(T)$ by
$$\|f\|^2_{T}=\|f\|^2+\|Tf\|^2, \ \ \ \ f\in \mathcal{D}(T).$$
A densely defined operator $T$ on $\mathcal{H}$ is said to be normal if $T$ is closed and $T^{\ast}T=TT^{\ast}$. \\
%
%
%
%
%
%
%

%

\section{ \sc\bf Unbounded weighted conditional expectation  }
Let $(X,\Sigma,\mu)$ be a $\sigma$-finite measure space and let
$\mathcal{A}$ be a $\sigma$-subalgebra of $\Sigma$ such that
$(X,\mathcal{A},\mu)$ is also $\sigma$-finite. We denote the
collection of (equivalence classes modulo sets of zero measure of)
$\Sigma$-measurable complex-valued functions on $X$ by
$L^0(\Sigma)$ and the support of a function $f\in L^0(\Sigma)$ is
defined as $S(f)=\{x\in X; f(x)\neq 0\}$. Moreover, we set
$L^2(\Sigma)=L^2(X,\Sigma,\mu)$. We also
adopt the convention that all comparisons between two functions or
two sets are to be interpreted as holding up to a $\mu$-null set.\\

For each $\sigma$-finite subalgebra $\mathcal{A}$ of $\Sigma$, the
conditional expectation, $E^{\mathcal{A}}(f)$, of $f$ with respect
to $\mathcal{A}$ is defined whenever $f\geq0$ almost everywhere or
$f\in L^2$. In any case, $E^{\mathcal{A}}(f)$
is the unique $\mathcal{A}$-measurable function for which
$$\int_{A}fd\mu=\int_{A}E^{\mathcal{A}}fd\mu, \ \ \ \forall A\in \mathcal{A} .$$
 As an operator on
$L^{2}({\Sigma})$, $E^{\mathcal{A}}$ is an idempotent and
$E^{\mathcal{A}}(L^2(\Sigma))=L^2(\mathcal{A})$. If there is no
possibility of confusion we write $E(f)$ in place of
$E^{\mathcal{A}}(f)$ \cite{rao,z}.\\

Let $u$ be a complex $\Sigma$-measurable function on $X$. Define
the measure $\mu_u:\Sigma\rightarrow [0, \infty]$  by

$$\mu_u(E)=\int_{E}|u|^2d\mu, \ \ \ \ E\in \Sigma.$$

It is clear that the measure $\mu_u$ is also $\sigma$-finite.\\


In this paper we consider $u$ is conditionable (i.e., $E(u)$ is defined). Operators of the form $EM_u(f)=E(u.f)$ acting in
$L^2(\mu)=L^2(X,\Sigma,\mu)$ with $\mathcal{D}(EM_u)=\{f\in
L^2(\mu):E(u.f) \in L^2(\mu)\}$ are called weighted
conditional expectation type operators. First we give some conditions under which the operator $EM_u$ is densely defined.\\

\vspace*{0.3cm} {\bf Lemma 2.1.} Let $w=1+E(|u|^2)$,  $\mu$ and $d\nu=wd\mu$. We get that\\

(i) $S(w)=X$ and $L^2(\nu)=L^2(X,\Sigma,\nu)\subseteq\mathcal{D}(EM_u)$.\\

(ii) If $E(|u|^2)<\infty$ a.e., then\\

$$\overline{L^2(\nu)}^{\|.\|_{\mu}}=\overline{\mathcal{D}(EM_u)}^{\|.\|_{\mu}}=L^2(\mu).$$\\

\vspace*{0.3cm} {\bf Proof.} Let $f\in L^2(\nu)$. Then
\begin{align*}
\|f\|^2_{\mu}&=\int_{X}|f|^2d\mu\\
&\leq \int_{X}|f|^2d\mu+\int_{X}E(|u|^2)|f|^2d\mu\\
&=\|f\|^2_{\nu}<\infty,\\
\end{align*}
so $f\in L^2(\mu)$. Also, by conditional-type H\"{o}lder-inequality we have
\begin{align*}
\|EM_u(f)\|^2_{\mu}&=\int_{X}|E(uf)|^2d\mu\\
&\leq\int_{X}E(|u|^2)E(|f|^2)d\mu\\
&=\int_{X}E(|u|^2)|f|^2d\mu\\
&\leq\|f\|^2_{\nu}<\infty,\\
\end{align*}
this implies that $f\in\mathcal{D}(EM_u)$.\\

Now we prove that $L^2(\nu)$ is dense in $L^2(\mu)$. Suppose that $f\in L^2(\mu)$ such that $\langle f,g \rangle=\int_{X}f.\bar{g}d\mu=0$ for all $g\in L^2(\nu)$. For $A\in \Sigma$ we set $A_n=\{x\in A:w(x)\leq n\}$. It is clear that $A_n\subseteq A_{n+1}$ and $X=\cup^{\infty}_{n=1}A_n$. Also, $X$ is $\sigma$-finite, hence $X=\cup^{\infty}_{n=1}X_n$ with $\mu(X_n)<\infty$. If we set $B_n=A_n\cap X_n$, then $B_n\nearrow A$ and so $f.\chi_{B_n}\nearrow f.\chi_{A}$ a.e. $\mu$. Since $\nu(B_n)\leq (n+1) \mu(B_n)<\infty$, we have $\chi_{B_n}\in L^2(\nu)$ and by our assumption $\int_{B_n}fd\mu=0$. Therefore by Fatou's lemma we get that $\int_{A}fd\mu=0$. Thus for all $A\in \Sigma$ we have $\int_{A}fd\mu=0$. This means that $f=0$ a.e. $\mu$ and so $L^2(\nu)$ is dense in $L^2(\mu)$.\\

\vspace*{0.3cm} {\bf Proposition 2.2.} If $u:X\rightarrow \mathbb{C}$ is $\Sigma$-measurable, then the following conditions are equivalent:\\

(i) $EM_u$ is densely defined,\\

(ii) $w-1=E(|u|^2)<\infty$ a.e. $\mu$,\\

(ii) $\mu_{E(|u|^2)}\mid_{\mathcal{A}}$ is $\sigma$-finite.\\

\vspace*{0.3cm} {\bf Proof.} $(i)\rightarrow (ii)$ Set $E=\{E(|u|^2)=\infty\}$. Clearly by Lemma 2.1, $f\mid_{E}=0$ a.e. $\mu$ for every $f\in L^2(\nu)$. This implies that  $f.w\mid_{E}=0$ a.e. $\mu$ for every $f\in L^2(\mu)$. So we have $w.\chi_{A\cap E}=0$ a.e. $\mu$ for all $A\in \Sigma$ with $\mu(A)<\infty$. By the $\sigma$-finiteness of $\mu$ we have $w.\chi_{E}=0$ a.e. $\mu$. Since $S(w)=X$, we get that $\mu(E)=0$.\\
$(ii)\rightarrow (i)$ It is clear by Lemma 2.1.\\

$(ii)\rightarrow (iii)$ Since $(X,\mathcal{A}, \mu)$ is a $\sigma$-finite measure space and $E(|u|^2)<\infty$ a.e. $\mu$, then there exists a sequence $\{X_n\}^{\infty}_{n=1}\subseteq \mathcal{A}$ such that $\mu(X_n)<\infty$ and $E(|u|^2)<n$ a.e. $\mu$ on $X_n$ for every $n\in \mathbb{N}$ and $X_n\nearrow X$ as $n\rightarrow \infty$.\\So we have
\begin{align*}
\mu_{(w-1)}\mid_{\mathcal{A}}(X_n)=\int_{X_n}E(|u|^2)d\mu\leq n\mu(X_n)<\infty, \ \ \ \ n\in \mathbb{N}.\\
\end{align*}
This yields (iii).\\

$(iii)\rightarrow (i)$ Let $\{A_n\}^{\infty}_{n=1}\subseteq \mathcal{A}$ be a sequence such that $A_n\nearrow X$ as $n\rightarrow \infty$ and $\mu_{w-1}\mid_{\mathcal{A}}(A_n)<\infty$ for every $k\in \mathbb{N}$. It follows from the definition of $\mu_{w-1}$ that $w-1=E(|u|^2)<\infty$ a.e. $\mu$ on $X$. Applying Lemma 2.1 we obtain $(i)$.\\

It is easily seen that $\mathcal{D}(M_{\bar{u}}E)\subseteq \mathcal{D}(EM_u)$. Also, by the previous proposition easily we get that: the operator $EM_u$ is densely defined if and only if the operator $M_{\bar{u}}E$ is densely defined. In this case we get that $\overline{\mathcal{D}(EM_{u})}^{\|.\|_2}=\overline{\mathcal{D}(M_{\bar{u}}E)}^{\|.\|_2}$.\\

\vspace*{0.3cm} {\bf Lemma 2.3}  Let the linear
transformation $T=EM_u$ be densely defined on $L^2(\Sigma)$, then $T^{\ast}=M_{\bar{u}}E$.\\

\vspace*{0.3cm} {\bf Proof} Let $f\in \mathcal{D}(T)$  and $g\in \mathcal{D}(T^{\ast})$. So we have
\begin{align*}
\langle Tf,g\rangle&=\int_{X}E(uf)\bar{g}d\mu\\
&=\int_{X}E(\bar{g})E(uf)d\mu\\
&=\int_{X}fuE(\bar{g})d\mu\\
&=\int_{X}f\overline{\bar{u}E(g)}d\mu\\
&=\langle f,M_{\bar{u}}Eg\rangle.\\
\end{align*}

From this relation we conclude that $M_{\bar{u}}E\subseteq (EM_u)^{\ast}$ and so $\mathcal{D}((EM_u)^{\ast})\subseteq \mathcal{D}(M_{\bar{u}}E)$.
Since the inner product function is continuous and $EM_u$ is densely defined, then $\mathcal{D}((EM_u)^{\ast})=\mathcal{D}(M_{\bar{u}}E)$. This completes the proof of the equality $(EM_u)^{\ast}=M_{\bar{u}}E$. \\

\vspace*{0.3cm} {\bf Proposition 2.4.} If $E(|u|^2)<\infty$ a.e. $\mu$. Then the linear transformation $EM_u:\mathcal{D}(EM_u)\rightarrow L^2(\Sigma)$ is closed.\\

\vspace*{0.3cm} {\bf Proof} Assume that $f_n\in \mathcal{D}(EM_{u})$, $f_n\rightarrow f$, $E(uf_n)\rightarrow g$, and let $h\in  \mathcal{D}(M_{\bar{u}}E)$. Then

\begin{align*}
\langle f,M_{\bar{u}}Eh\rangle&=\lim_{n\rightarrow \infty}\langle f_n,M_{\bar{u}}Eh\rangle\\
&=\lim_{n\rightarrow \infty}\langle E(uf_n),h\rangle\\
&=\langle g,h\rangle.\\
\end{align*}

This calculation ( which uses the continuity of the inner product and the fact that $f_n\in \mathcal{D}(EM_{u})$) shows that $f\in \mathcal{D}(EM_{u})$ and $E(uf)=g$, as required.\\


\vspace*{0.3cm} {\bf Proposition 2.5.} If $\mathcal{D}(EM_u)$
is dense in $L^2(\Sigma)$ and $u$ is almost every where finite valued, then the operator $EM_u$ is normal if and only if $u\in L^0(\mathcal{A})$.\\

\vspace*{0.3cm} {\bf Proof.} Let the operator $EM_u$ be normal. Then for every $f\in \mathcal{D}(M_{\bar{u}}EM_u)=\mathcal{D}(EM_{E(|u|^2)})$ we have
$$T^{\ast}Tf=TT^{\ast}f \  \ \mu, \ a.e., \ \ \Rightarrow \ \ \bar{u}E(uf)=E(|u|^2)E(f) \ \ \mu, \ a.e.,$$

by taking $E$ over both side of the equality, for a positive element $a$ of $L^2(X,\mathcal{A},\nu)$ we get that $|E(u.\chi_{A_n})|^2a=E(|u|^2.\chi_{A_n})a \ \ \ \mu, \ a.e.,$ in which $A_n=\{x\in X:|u(x)|\leq n\}$ and $\mu(A_n)<\infty$. This
implies that $|E(u.\chi_{A_n})|^2=E(|u|^2.\chi_{A_n})$. Since $E(|u.\chi_{A_n}-E(u.\chi_{A_n})|^2)=E(|u.\chi_{A_n}|^2)-|E(u.\chi_{A_n})|^2$, then $u.\chi_{A_n}=E(u.\chi_{A_n})$ and so $u=E(u)$, i.e, $u\in L^0(\mathcal{A})$.\\

Conversely. If $EM_u$ is densely defined, then by the same method of Lemma 2.1 we get that

$$L^2(\nu)\subseteq \mathcal{D}(M_{\bar{u}}EM_u), \ \ \ \ \ L^2(\nu)\subseteq \mathcal{D}(EM_{E(|u|^2)}), \ \ \ \mathcal{D}(EM_{E(|u|^2)})\subseteq \mathcal{D}(M_{\bar{u}}EM_u)$$
and
$$\overline{\mathcal{D}(M_{\bar{u}}EM_u)}=\overline{L^2(\nu)}=\overline{\mathcal{D}(EM_{E(|u|^2)})}=L^2(\mu),$$
where $d\nu=(1+(E(|u|^2))^2)d\mu$. Suppose that $u\in L^0(\mathcal{A})$. Then $\mathcal{D}(M_{\bar{u}}EM_u)=\mathcal{D}(EM_{E(|u|^2)})=\mathcal{D}(M_{|u|^2}E)$ and for every $f\in \mathcal{D}(M_{|u|^2}E)$ we have

$$T^{\ast}Tf=|u|^2E(f)=TT^{\ast}f \  \ \mu, \ a.e.$$

On the other hand, since $EM_u$ is densely defined. Then by Lemma 2.4 it is closed. This implies that the operator $EM_u$ is normal.\\

\vspace*{0.3cm} {\bf Proposition 2.6.} If $\mathcal{D}(EM_u)$
is dense in $L^2(\Sigma)$, then the operator $EM_u$ is self-adjoint if and only if $u\in L^0(\mathcal{A})$ is real valued.\\

\vspace*{0.3cm} {\bf Proof.} Suppose that the operator $EM_u$ is self-adjoint, then it is normal. So by Proposition 2.5 we get that $u\in L^0(\mathcal{A})$. Let $a$ be a positive $\mathcal{A}$-measurable function in $\mathcal{D}(EM_u)$. Hence
$$\bar{u}a=\bar{u}E(a)=T^{\ast}(a)=T(a)=E(ua)=ua,$$
this implies that $\bar{u}=u$ i.e, $u$ is real valued.\\
Conversely, if $u\in L^0(\mathcal{A})$ is real valued, then for $f, g\in \mathcal{D}(M_{\bar{u}}E)\subseteq \mathcal{D}(EM_u)$,
\begin{align*}
\langle E(uf),g\rangle&=\int_{X}E(uf)\bar{g}d\mu\\
&=\int_{X}uE(f)\bar{g}d\mu\\
&=\int_{X}\bar{u}E(f)\bar{g}d\mu\\
&=\langle M_{\bar{u}}E(f),g\rangle.\
\end{align*}
This implies that $\mathcal{D}(EM_u)=\mathcal{D}(M_{\bar{u}}E)$ and $T^{\ast}(f)=T(f)$ for all $f\in \mathcal{D}(EM_u)$.\\

It is well-known that for a densely defined closed operator $T$ of $\mathcal{H}_1$ into $\mathcal{H}_2$, there exists a partial isometry $U_T$ with initial space $\mathcal{N}(T)^{\perp}=\overline{\mathcal{R}(T^{\ast})}=\overline{\mathcal{R}(|T|)}$
and final space $\mathcal{N}(T^{\ast})^{\perp}=\overline{\mathcal{R}(T)}$ such that $$T=U_T|T|.$$

A closed densely defined operator $T$ in $\mathcal{H}$ is said to be quasinormal if $U|T|\subseteq|T|U$, where $T=U|T|$ is the polar decomposition of $T$ \cite{br,js}. As is shown in [\cite{jjs}, Theorem 3.1], A closed densely defined operator $T$ in $\mathcal{H}$ is quasinormal if and only if $T|T|^2=|T|^2T$.\\

\vspace*{0.3cm} {\bf Proposition 2.7.} Suppose that  $\mathcal{D}(EM_u)$
is dense in $L^2(\Sigma)$. Let $EM_u=U|EM_u|$ be the polar decomposition of $EM_u$. Then\\

(i) $|EM_u|=M_{u'}EM_u$, where $u'=(E(|u|^2))^{\frac{-1}{2}}.\chi_{S}.\bar{u}$ and $S=S(E(|u|^2))$,\\

(ii) $U=EM_{\widetilde{u}}$, where $\widetilde{u}:X\rightarrow \mathbb{C}$ is an a.e. $\mu$ well-defined $\Sigma$-measurable function such that

$$\widetilde{u}=u.\frac{1}{(E(|u|^2))^{\frac{1}{2}}}.\chi_{S}.$$

\vspace*{0.3cm} {\bf Proof.} (i). For every $f\in \mathcal{D}(M_{u'}EM_u)$ we have
\begin{align*}
\|M_{u'}EM_u(f)(f)\|^2&=\int_{X}((|u|^2))^{-1}.\chi_{S}|u|^2|E(uf)|^2d\mu\\
&=\int_{X}((|u|^2))^{-1}.\chi_{S}E(|u|^2)|E(uf)|^2d\mu\\
&=\||EM_u|(f)\|^2.\\
\end{align*}
Also, by Lemma 2.1 we conclude that $\mathcal{D}(M_{u'}EM_u)=\mathcal{D}(|EM_u|)$ and it is easily seen that $M_{u'}EM_u$ is a positive operator. This observations imply that $|EM_u|=M_{u'}EM_u$.\\

(ii). For $f\in L^2(\Sigma)$ we have
$$\int_{X}|E(\widetilde{u}f)|^2d\mu=\int_{X}\frac{1}{(E(|u|^2))}.\chi_{S}|E(uf)|^2d\mu,$$
which implies that the operator $EM_{\widetilde{u}}$ is well-defined and $\mathcal{N}(EM_u)=\mathcal{N}(EM_{\widetilde{u}})$. Also, for $f\in \mathcal{D}(EM_u)\ominus \mathcal{N}(EM_u)$ we have
\begin{align*}
U(|EM_u|(f))&=\frac{1}{(E(|u|^2))^{\frac{1}{2}}}.\chi_{S}E(u(E(|u|^2))^{\frac{-1}{2}}\chi_{S}\bar{u}E(uf))\\
&=\frac{1}{E(|u|^2)}.\chi_{S}E(|u|^2)E(uf)\\
&=E(uf).\\
\end{align*}
Thus $\|U(f)\|=\|f\|$ for all $f\in \mathcal{R}(|EM_u|)$ and since $U$ is a contraction, then it holds for all $f\in \mathcal{N}(EM_u)^{\perp}=\overline{\mathcal{R}(|EM_u|)}$.\\

Here we remind that: A complex number $\lambda$ belongs  to the resolvent set $\rho(T)$ of the closed linear operator $T$ on a Hilbert space $\mathcal{H}$, if the operator $\lambda I-T$ has a bounded everywhere on $\mathcal{H}$ defined inverse $(\lambda I-T)^{-1}$, called the resolvent of $T$ at $\lambda$ and denoted by $R_{\lambda}(T)$.\\
The set $\sigma(T):=\mathbb{C}\setminus \rho(T)$ is called the spectrum of the operator $T$.\\

{\bf Theorem 2.8.} Let $EM_u$ be densely defined and $\mathcal{A}\varsubsetneq\Sigma$, then\\

(i) $\operatorname{ess range(E(u))}\cup\{0\}\subseteq \sigma(EM_u)$,\\

(ii) If $L^2(\mathcal{A})\subseteq \mathcal{D}(EM_u)$, then $\sigma(EM_u)\subseteq \operatorname{ess range(E(u))}\cup\{0\}$.\\

{\bf Proof.} (i) Since $\mathcal{A}\neq\Sigma$, we get that $EM_u$ is not surjective. So $0\in \sigma(EM_u)$. Let $0\neq\lambda\in \operatorname{ess range(E(u))}$. Then the measure of the set

$$N_{n}:=\{x\in X: |E(u)(x)-\lambda|<n^{-1}\}$$
is positive for each $n\in \mathbb{N}$. Because $\mu$ is $\sigma$-finite, upon replacing $N_{n}$  by a subset, we can assume that $0<\mu(N_n)<\infty$.  Since $E(u)$ is bounded on $N_n$, $\chi_{N_n}\in \mathcal{D}(EM_u)$ and
\begin{align*}
\|(\lambda I-EM_u)\chi_{N_{n}}\|^2&=\int_{X}|\lambda \chi_{N_n}-E(u\chi_{N_n})|^2d\mu\\
&=\int_{X}|\lambda \chi_{N_n}-E(u)\chi_{N_n}|^2d\mu\\
&=\int_{N_n}|(\lambda -E(u))|^2d\mu\\
&\leq n^{-2}\|\chi_{N_{n}}\|^2.
\end{align*}
Hence $\lambda$ is not a regular point. Therefore, $\lambda\notin \rho(EM_u)$, and so $\lambda\in \sigma(EM_u)$.\\

(ii) Suppose that
$\lambda\notin \operatorname{ess range(E(u))}$. Then there exists $n\in \mathbb{N}$ such that $\mu(N_n)=0$. We
show that $T-\lambda I$ is invertible. If $Tf-\lambda f=0$ for $f \in \mathcal{D}(EM_u)$, then
$E(uf)=\lambda f$. So $f$ is $\mathcal{A}-$measurable. Thus
$(E(u)-\lambda)f=E(uf)-\lambda f=0$. Since $\lambda\notin
\operatorname{ess range(E(u))}$, then
$E(u)-\lambda\geq\varepsilon$ a.e for some $\varepsilon>0$. So
$f=0$ a.e. This implies that $T-\lambda I$ is injective and so is invertible.\\
 Now we
show that $T-\lambda I$ is surjective. Let $g\in
L^{2}(\Sigma)$. We can write

$$g=g-E(g)+E(g), \ \ \ g_{1}=g-E(g),\ \ \ g_{2}=E(g).$$

Clearly we have $g_{2}\in L^2(\mathcal{A})\subseteq \mathcal{D}(EM_u)$ and  $g_{1}\in L^{2}(\Sigma)$,
$(E(g_{1})=0$. Let $$f_{1}=\frac{\lambda
g_{1}+T(g_{2})}{\lambda(E(u)-\lambda)}, \ \ \
f_{2}=\frac{-g_{2}}{\lambda}.$$ Since $\lambda\notin
\operatorname{ess range(E(u))}$, we get
$E(u)-\lambda\geq\varepsilon$ a.e for some $\varepsilon>0$. So
$\|\frac{1}{E(u)-\lambda}\|_{\infty}\leq\frac{1}{\varepsilon}$.
Thus $f_{2}\in L^{2}(\mathcal{A})$,  $f_{1}\in
L^{2}(\Sigma)$ and $f=f_{1}+f_{2}\in  L^{2}(\Sigma)$. Direct
computation shows that $T(f)-\lambda f=g$.
This implies that
$T-\lambda I$ is invertible and so $\lambda \notin \sigma(T)$.\\

\vspace*{0.3cm} {\bf  Remark.}
If $EM_u$ is every where defined, then $\operatorname{ess range(E(u))}\cup\{0\}=\sigma(EM_u)$.\\



\vspace*{0.3cm} {\bf Proposition 2.9.} Let $EM_u$ be densely defined on $L^2(\Sigma)$. Then we have\\

 (i) If $\bar{u}E(u)=E(|u|^2)$, $\mu$, a.e,. on $S(E(|u|^2))$, then the operator $EM_u$ is quasinormal.\\

 (ii) If the operator $EM_u$ is quasinormal, then $\bar{u}E(u)=E(|u|^2)$, $\mu$, a.e,. on $S(E(u))$.\\

\vspace*{0.3cm} {\bf proof.} (i) Let $\nu_1=1+(E(|u|^2))^3$ and $\nu_2=1+|E(u)|^2(E(|u|^2))^2$. It is obvious that $\nu_2\leq \nu_1$, $\mu$, a.e,. So $L^2(\nu_1)\subseteq L^2(\nu_2)\subseteq L^2(\mu)$. If we set  $T_1=EM_u|EM_u|^2=EM_{uE(|u|^2)}$ and $T_2=|EM_u|^2EM_u=M_{\bar{u}E(u)}EM_u$, then by Proposition 2.2 we conclude that $EM_u$ is densely defined if and only if $T_1$ is densely defined. Also, it is easily seen that $\mathcal{D}(T_1)\subseteq \mathcal{D}(T_2)$. So,
$$\overline{L^2(\nu_1)}=\overline{L^2(\nu_2)}=\overline{\mathcal{D}(T_1)}=\overline{\mathcal{D}(T_2)}=L^2(\mu).$$
Since $\bar{u}E(u)=E(|u|^2)$, $\mu$, a.e,. on $S(E(|u|^2))$, then $|E(u)|^2=E(|u|^2)$. So for all $f\in \mathcal{D}(T_2)$ we get that
\begin{align*}
\|T_1(f)\|^2&=\int_{X}(E(|u|^2))^2|E(uf)|^2d\mu\\
&=\int_{X}|E(u)|^2E(|u|^2)|E(uf)|^2d\mu\\
&=\|T_2(f)\|^2.\\
\end{align*}
This implies that $\mathcal{D}(T_1)=\mathcal{D}(T_2)$. Also, for all $f\in \mathcal{D}(T_2)$ we have
\begin{align*}
T_2&=M_{\bar{u}E(u)}EM_u(f)\\
&=\bar{u}E(u)E(uf)\\
&=E(|u|^2)E(uf)\\
&=T_1(f).\\
\end{align*}

(ii) Suppose that the operator $EM_u$ is quasinormal, then for all $f\in\mathcal{D}(|EM_u|^2EM_u)$ we have
\begin{align*}
EM_u|EM_u|^2(f)&=E(uE(|u|^2)f)\\
&=E(|u|^2)E(uf)\\
&=\bar{u}E(u)E(uf)\\
&=|EM_u|^2EM_u,\\
\end{align*}

So, for a positive $\mathcal{A}$-measurable function $a\in L^2(\nu_1)$ we get that
$\bar{u}E(u)E(u)a=E(|u|^2)E(u)a$ and then $\bar{u}E(u)=E(|u|^2)$, $\mu$, a.e,. on $S(E(u))$.\\

\vspace*{0.3cm} {\bf Corollary 2.10.} Let $EM_u$ be densely defined on $L^2(\Sigma)$ and $S(E(u))=S(E(|u|^2))$. Then the operator $EM_u$ is quasinormal if and only if $\bar{u}E(u)=E(|u|^2)$, $\mu$, a.e,.\\

The next proposition can be easily deduced from the closed graph
theorem.\\

\vspace*{0.3cm} {\bf Proposition 2.11.} If $T$ is a closed operator
on $\mathcal{H}$ such that $T(\mathcal{D}(T)\subseteq
(\mathcal{D}(T)$, then $T$ is a bounded operator on the Hilbert
space $(\mathcal{D}(T), \|.\|_{T})$.\\

\vspace*{0.3cm} {\bf Proposition 2.12.} Let $\mathcal{D}(EM_u)$
is dense in $L^2(\Sigma)$. Then the following conditions are
valid:\\

(i) If $EM_u(\mathcal{D}(EM_u))\subseteq
\mathcal{D}(EM_u)$, then $|E(u)|^4\leq c(1+|E(u)|^2)$ a.e. $\mu$.\\

(ii) Assume that, there exists $c>0$ such that $|E(u)|^2E(|u|^2)\leq c(1+|E(u)|^2)$ a.e. $\mu$ and $\int_{X}|E(u)f|^2d\mu\leq \int_{X}|E(uf)|^2d\mu$, for all $f\in \mathcal{D}(EM_u)$. Then $EM_u(\mathcal{D}(EM_u))\subseteq
\mathcal{D}(EM_u)$.\\

\vspace*{0.3cm} {\bf Proof.}$(i)$. Since
$EM_u$ is closed, densely defined and
$EM_u(\mathcal{D}(EM_u))\subseteq
\mathcal{D}(EM_u)$, then by closed graph theorem $EM_u$
is a bounded operator on
$(\mathcal{D}(EM_u),\|.\|_{EM_u})$.  Hence there exists
$c>0$ such that $\|EM_u(f)\|^2_{EM_u}\leq
c\|f\|^2_{EM_u}$ for $f\in \mathcal{D}(EM_u)$. By
replacing $f$ with $EM_u(f)$ we have
\begin{align*}
\|(EM_u)^2(f)\|^2&\leq \|EM_u(f)\|^2+\|(EM_u)^2(f)\|^2\\
&\leq c(\|f\|^2+\|EM_u(f)\|^2),\\
\end{align*}

i.e,
\begin{align*}
\int_{X}|E(u)|^2|E(uf)|^2d\mu\leq c(\int_{X}|f|^2d\mu+\int_{X}|E(uf)|^2d\mu).\\
\end{align*}

Thus for every $A\in \mathcal{A}$ with $\chi_{A}\in \mathcal{D}(EM_u)$  we have
\begin{align*}
\int_{A}|E(u)|^2|E(u)|^2d\mu\leq c(\int_{A}d\mu+\int_{A}|E(u)|^2d\mu).\\
\end{align*}

Since $EM_u$ is densely-defined, we get that $|E(u)|^4\leq c(1+|E(u)|^2)$.\\

$(ii)$. Let $f\in\mathcal{D}(EM_u)$. Then by
assumptions $|E(u)|^2E(|u|^2)\leq c(1+|E(u)|^2)$ a.e. $\mu$ and $\int_{X}|E(u)f|^2d\mu\leq \int_{X}|E(uf)|^2d\mu$, a.e. $\mu$, we have
\begin{align*}
\int_{X}|(EM_u)^2(f)|^2d\mu&=\int_{X}|E(u)|^2|E(uf)|^2d\mu\\
&\leq \int_{X}|E(u)|^2E(|u|^2)|f|^2d\mu\\
&\leq
c(\int_{X}|f|^2d\mu
+\int_{X}|E(u)f|^2d\mu)\\
&\leq c(\int_{X}|f|^2d\mu
+\int_{X}|E(uf)|^2d\mu)\\
& \leq c(\|f\|^2+\|E(uf)\|^2)\\
&<\infty.\\
\end{align*}
Therefore $EM_u(f)\in \mathcal{D}(EM_u)$.\\


\vspace*{0.3cm} {\bf Corollary 2.13.} If $\mathcal{D}(M_u)$ is dense in
$L^2(\Sigma)$, then the following conditions are
equivalent:\\

(i) $M_u(\mathcal{D}(M_u))\subseteq
\mathcal{D}(M_u)$.\\

(ii) There exists $c>0$ such that $u^4\leq c(1+u^2)$ a.e. $\mu$.\\

Here we present some examples of conditional
expectations and corresponding multiplication operators to
illustrate concrete application of the main results of the paper in this section.\\

 {\bf Example (i).} In this example we consider some cases of sub-algebras of the $\sigma$-algebra $\Sigma$ when $\mu(X)<\infty$. In this cases, we restrict our attention to weighted conditional type operators on sub-spaces of $L^2(\Sigma)$.\\

 Case 1. As our first cases, take the extreme case when $\mathcal{A}=\Sigma$. Then $E(f)=f$ for all $f\in L^2(\Sigma)$; that is, $E=I$, where $I$ is the identity operator on $L^2(\Sigma)$. In this situation the weighted conditional type operators are simply multiplication operators: $EM_u=M_u$. Therefore our results contain the similar standard results regarding multiplication operators as follows:\\

 (a) $M_u$ is densely defined if and only if $u<\infty$ a.e., with respect to $\mu$,\\

 (b) If $u<\infty$ a.e., $\mu$, then\\

 $(b)_1$ $M_u$ is self adjoint if and only if $u$ is real valued,\\

 $(b)_2$ $M_u$ is closed,\\

 $(b)_3$  $\operatorname{ess range(u)}=\sigma(M_u)$,\\

  etc.\\

 Case 2. If $\mathcal{A}=\{X, \emptyset\}$. Here $\mathcal{A}$-measurable functions are constant on $X$ and
 $$E(f)=(\mu(X))^{-1}\int_{X}fd\mu, \ \ \ \ \ f\in L^2(\Sigma).$$

  In this case the weighted conditional type operator $EM_u$ acts as

  $$EM_u(f)=(\mu(X))^{-1}\int_{X}ufd\mu, \ \ \ \ \ f\in \mathcal{D}(EM_u).$$
   Again, we have:\\

 (a) $EM_u$ is densely defined if and only if $E(|u|^2)=(\mu(X))^{-1}\int_{X}|u|^2d\mu<\infty$ a.e., $\mu$, ( equivalently; if and only if $u$ is an $L^2$ function,\\

 (b) If $u$ is an $L^2$ function, then\\

 $(b)_1$ $EM_u$ is normal  if and only if $u=(\mu(X))^{-1}\int_{X}ud\mu$ a.e., $\mu$,\\

 $(b)_2$ $EM_u$ is self adjoint if and only if $u=(\mu(X))^{-1}\int_{X}ud\mu$ a.e., $\mu$ and $(\mu(X))^{-1}\int_{X}ud\mu$ is real,\\

 $(b)_3$ $EM_u$ is closed,\\

 $(b)_4$ $\sigma(EM_u)=\{(\mu(X))^{-1}\int_{X}ud\mu\}$,\\

  etc.\\

  Case 3. Suppose that the $\sigma$-sub-algebra $\mathcal{A}$ is generated by a countable partition $\{A_n:n\in \mathbb{N}\}$ of $X$ into disjoint sets of finite measure. It is known that the conditional expectation of any $f\in L^2(\Sigma)$ relative to $\mathcal{A}$ is:
  $$E(f)=\sum^{\infty}_{n=1}(\mu(A_n))^{-1}(\int_{A_n}fd\mu).\chi_{A_n}.$$
In this case the weighted conditional type operator $EM_u$ acts as

  $$EM_u(f)=\sum^{\infty}_{n=1}(\mu(A_n))^{-1}(\int_{A_n}ufd\mu).\chi_{A_n}, \ \ \ \ \ f\in \mathcal{D}(EM_u).$$

Then we have:\\

 (a) $EM_u$ is densely defined if and only if $E(|u|^2)=\sum^{\infty}_{n=1}\beta_n \chi_{A_n}<\infty$ a.e., $\mu$, where $\beta_n=(\mu(A_n))^{-1}(\int_{A_n}|u|^2d\mu)$,\\

 (b) If $\sum^{\infty}_{n=1}\beta_n \chi_{A_n}<\infty$ a.e., $\mu$, then\\

 $(b)_1$ $EM_u$ is normal  if and only if $u=\sum^{\infty}_{n=1}(\mu(A_n))^{-1}(\int_{A_n}ud\mu).\chi_{A_n}$ a.e., $\mu$,\\

 $(b)_2$ $EM_u$ is self adjoint if and only if $u=\sum^{\infty}_{n=1}(\mu(A_n))^{-1}(\int_{A_n}ud\mu).\chi_{A_n}$ a.e., $\mu$ and $u$ is real,\\

 $(b)_3$ $EM_u$ is closed,\\

 $(b)_4$ $\sigma(EM_u)=\{\beta_n:n\in \mathbb{N}\}$,\\

  etc.\\

{\bf Example (ii).} Let $X=[0,1]\times [0,1]$, $d\mu=dxdy$,
$\Sigma$ the Lebesgue subsets of $X$ and let
$\mathcal{A}=\{A\times [0,1]: A \ \mbox{is a Lebesgue set in} \
[0,1]\}$. Then, for each $f$ in $L^2(\Sigma)$, $(Ef)(x,
y)=\int_0^1f(x,t)dt$, which is independent of the second
coordinate. Then we have:\\

 (a) $EM_u$ is densely defined if and only if $E(|u|^2)(x,y)=\int_0^1|u(x,t)|^2dt<\infty$ a.e., $\mu$,\\

 (b) If $\int_0^1u(x,t)dt<\infty$ a.e., $\mu$, then\\

 $(b)_1$ $EM_u$ is normal  if and only if $u(x,y)=\int_0^1u(x,t)dt$ a.e., $\mu$,\\

 $(b)_2$ $EM_u$ is self adjoint if and only if $u(x,y)=\int_0^1u(x,t)dt$ a.e., $\mu$ and $u$ is real valued,\\

 $(b)_3$ $EM_u$ is closed,\\

 $(b)_4$ $\sigma(EM_u)=\{\int_0^1u(x,t)dt:x\in E\}$, where $E=\{0\leq x\leq1:\int_0^1u(x,t)dt<\infty\}$,\\

  etc.\\

{\bf Example (iii).} Let $\Omega=[-1,1]$, $d\mu=\frac{1}{2}dx$ and $\mathcal{A}=<\{(-a,a):0\leq a\leq1\}>$ (Sigma algebra generated by symmetric intervals).
Then
 $$E^{\mathcal{A}}(f)(x)=\frac{f(x)+f(-x)}{2}, \ \ x\in \Omega,$$
 where $E^{\mathcal{A}}(f)$ is defined. If we set $u(x)=e^x$, then $E(|u|^2)(x)=\cosh(2x)$ and so\\

 (a) $EM_u$ is densely defined,\\

 (b) $EM_u$ can not be normal at all, since $e^x\neq\cosh(x)$ for every $x\neq 0$,\\

 (c) $EM_u$ can not be self adjoint at all, since $e^x\neq\cosh(x)$ for every $x\neq 0$,\\

 (d) $EM_u$ is closed,\\

 (e) $\sigma(EM_u)=\{\cosh(x):-1\leq x\leq 1\}$,\\

  etc.\\

{\bf Example (iv).}  Let $X=\mathbb{N}\cup\{0\}$,
$\mathcal{G}=2^{\mathbb{N}}$ and let
$\mu(\{x\})=\frac{e^{-\theta}\theta^x}{x !}$, for each $x\in X$
and $\theta\geq0$. Elementary calculations show that $\mu$ is a
probability measure on $\mathcal{G}$. Let $\mathcal{A}$ be the
$\sigma$-algebra generated by the partition $B=\{\emptyset, X,
\{0\}, X_1=\{1, 3, 5, 7, 9, ....\}, X_2=\{2, 4, 6, 8, ....\},\}$
of $\mathbb{N}$. Note that $\mathcal{A}$ is a sub-$\sigma$-finite
algebra of $\Sigma$ and each of element of $\mathcal{A}$ is an
$\mathcal{A}$-atom. Thus the conditional expectation of any $f\in
\mathcal{D}(E)$ relative to $\mathcal{A}$ is constant on
$\mathcal{A}$-atoms. Hence there exists scalars $a_1, a_2, a_3$
such that
$$E(f)=a_1\chi_{\{0\}}+a_2\chi_{X_1}+a_3\chi_{X_2}.$$
So
$$E(f)(0)=a_1, \ \ \ \ E(f)(2n-1)=a_2, \ \ \ \ E(f)(2n)=a_3,$$
for all $n\in \mathbb{N}$. By definition of conditional
expectation with respect to $\mathcal{A}$, we have

 $$a_1=f(0), \ \ \ \ a_2=\frac{\sum_{n\in \mathbb{N}}f(2n-1)\frac{e^{-\theta}\theta^{2n-1}}{(2n-1)
!}}{\sum_{n\in \mathbb{N}}\frac{e^{-\theta}\theta^{2n-1}}{(2n-1)
!}}, \ \ \ \ a_3=\frac{\sum_{n\in
\mathbb{N}}f(2n)\frac{e^{-\theta}\theta^{2n}}{(2n) !}}{\sum_{n\in
\mathbb{N}}\frac{e^{-\theta}\theta^{2n}}{(2n) !}}.$$

For example, if we set $u(x)=x$, then $E(u)$ is a special function
as follows;
$$E(u)=\theta
\coth(\theta)\chi_{X_1}+\frac{\cosh(\theta)-1}{\cosh(\theta)}\chi_{X_2}.$$\\
And so we have:\\

(a) $EM_u$ is densely defined,\\

 (b) $EM_u$ is normal if and only if $EM_u$ is self adjoint if and only if
$ \theta=\theta \coth(\theta)$ when $\theta \in X_1$ and $\theta=\frac{\cosh(\theta)-1}{\cosh(\theta)}$ when $\theta \in X_2$.\\

(c) $EM_u$ is closed,\\

 (d) $\sigma(EM_u)=\{\theta
\coth(\theta),\frac{\cosh(\theta)-1}{\cosh(\theta)}\}$,\\

  etc.\\

\end{document}